\documentclass{article}
\usepackage{amssymb,amsmath}
\newtheorem{propo}{Proposition}[section]
\newtheorem{lema}[propo]{Lemma}
\newtheorem{definicion}[propo]{Definition}

\newtheorem{notacion}[propo]{Notation}

\newtheorem{teorema}[propo]{Theorem}
\newenvironment{prueba}[1]{\paragraph{\sl {\bf Proof}#1}}{\qed}

\def\un{\underline}
\def\ov{\overline}

\def\un{\underline}

\def\kk{{\bf k}}

\newcommand{\CC}{{\bf C}}
\newcommand{\RR}{{\bf R}}

\newcommand{\QQ}{{\bf Q}}
\newcommand{\ZZ}{{\bf Z}}
\newcommand{\NN}{{\bf N}}

\newcommand{\cA}{{\mathcal A}}
\newcommand{\cB}{{\cal B}}
\newcommand{\cD}{{\cal D}}
\newcommand{\cE}{{\cal E}}
\newcommand{\cF}{{\cal F}}
\newcommand{\cG}{{\cal G}}

\newcommand{\cK}{{\cal K}}
\newcommand{\cL}{{\cal L}}
\newcommand{\cM}{{\cal M}}
\newcommand{\cN}{{\cal N}}
\newcommand{\cO}{{\cal O}}

\newcommand{\cR}{{\cal R}}

\newcommand{\cV}{{\cal V}}
\newcommand{\cU}{{\cal U}}

\newcommand{\wL}{{\widetilde L}}

\newcommand{\qed}{{$\square$}}

\newcommand{\Ker}{{\rm Ker}}

\newcommand{\ord}{{\rm ord}}

\newcommand{\gr}{{\rm gr}}
\newcommand{\too}{\rightarrow}
\newcommand{\weight}{{\rm weight}}

\author{F.J. Castro-Jim\'enez and M. Granger}
\title{A flatness property for filtered $\mathcal D$-modules}

\date{}
\begin{document}
\maketitle

\begin{abstract} Let $\cM$ be a coherent module over the ring $\cD_X$ of linear differential operators on an
analytic manifold $X$ and let  $Z_1,\cdots , Z_k$ be $k$ germs of transverse hypersurfaces at a point $x\in X$.  The Malgrange-Kashiwara V-filtrations along these hypersurfaces, associated with a given presentation of  the germ of $\cM$
at $x$, give rise to a multifiltration $U_\bullet(\cM)$  of $\cM_x$ as in Sabbah's paper   \cite{Sabbah}
and to an analytic standard fan   in a way similar to
\cite{assi-castro-granger-jpaa-2}. We prove here that this standard
fan is adapted to the multifiltration, in the sense of C. Sabbah. This
result completes the proof of the existence of an adapted fan in
\cite{Sabbah}, for which the use of   \cite{Sab-Cas} is not possible.
\end{abstract}

\section{Introduction} Let us consider a coherent module $\cM$ over the ring $\cD_X$ of differential operators on an
analytic manifold $X$. For any smooth hypersurface $Z$ of $X$,
Malgrange and Kashiwara defined a filtration along $Z$ for $\cD_X$
and the notion of a good filtration for $\cM$. Given a set of
transverse smooth hypersurfaces $Z_1, \ldots, Z_k$, Sabbah
considered in \cite{Sabbah} multifiltrations of $\cM$ indexed by
$k$-uples of relative integers. To be precise he dealt with
linear combinations over ${\bf Q}^+$ of the filtrations $V^{(j)}$
along each hypersurface $Z_j$ and with refinements $V^{\Gamma}$ of
the original multifiltration associated with each rational
polyedral simplicial cone in the positive quadrant of $({\bf
Q}^k)^{\star}$. The original multifiltration is the one which
corresponds to the case $\Gamma ={\bf N}^k$.

The aim of this paper is to clarify the flatness properties which
appear in \cite{Sabbah}, namely to prove the existence of a fan
$\cE$, such that for any cone $\Gamma $ in this fan, the Rees
module of the filtration $V^{\Gamma}$ is flat over the affine ring
$A_{\Gamma}$ of the toric blowing-up of ${\CC}^k$ associated with
the fan $\cE$. Such a fan is called an adapted fan in
\cite{Sabbah}. The reason for this clarification is that  the proof in  \cite{Sabbah}
depends on  the appendix \cite{Sab-Cas} in which the main tool is
a  division theorem in Rees rings of differential operators which is not correct
as stated. Indeed, the infinite process that its proof suggests, would require
monomials of unbounded degrees in the differential variables of the ring of operators.

One of the main consequences of the existence of an adapted fan as
developed in \cite{Sabbah} is that we thus complete the proof of
the existence of non trivial functional equations of
Bernstein-Sato type for a $k$--uple of functions, following the
argument of Sabbah in \cite{Sabbah}. It should be emphasized here
that this problem involves a  $\cD _{X\times \CC ^k}$-module
naturally associated with a $k$-uple of functions on $X$, and that
in this situation a multifiltration, along the transverse
hypersurfaces $t_j=0$ due to the factor $\CC ^k$, appears in a
natural way.

We must note here that the proof of the existence of
Bernstein-Sato equations has already been completed in Bahloul's
paper \cite{bahloul}, by a different method which avoids the
reference to a flatness property involved in the notion of an
adapted fan. Bahloul uses instead the analytic standard fan as
defined in \cite{assi-castro-granger-jpaa-2}, and described also
in \cite{bahloul-takayama}. It is therefore not completely
surprising that the adapted fan wanted in \cite{Sabbah} turns out in
fact to be the analytic standard  fan. This emphasizes the
interest of Bahloul's proof which also has the advantage of being
constructive. This constructiveness is valuable in an
algebraic setting also. In both proofs, Bahloul's  and
Sabbah's , the latter which is completed by this paper, the main step
is the proof
of the goodness of the so-called saturated filtration. The details
are given in section \ref{bfunctions}.

The starting point of the proof of our main theorem (Theorem
\ref{mainresult}) is a criterion of flatness by M. Herrmann and U.
Orbanz (see \cite{Herrmann-Orbanz}), for graded modules over
graded rings, where the grading is indexed by an arbitrary
commutative group. In the statement of this result there is no
reference to any finiteness property. In our case the indexing
group will be $\ZZ ^k$, and the ring a conical subring of the ring
of Laurent polynomials in $k$ variables with the obvious
multigrading by monomials. The main ingredient of the proof of
Theorem \ref{mainresult} is then the existence of a simultaneous
$L$--standard basis of a submodule $\cN\subset \cD^r$ for  all $L$
in a cone of the analytic standard  fan of $\cN$ (see
\cite{assi-castro-granger-jpaa-2}).

The plan of the paper is as follows. In section \ref{multif_rees},
we recall with more details what $V$-multifiltrations, and
Rees modules are, and the refinements of these multifiltrations with
respect to a rational polyhedral simplicial cone. We end this
section by defining the fiber at zero of these Rees modules seen
as modules  over the ring of an affine toric blowing up. We aim to
prove that this is the fiber of a flat deformation if the cone is
taken in the analytic standard fan.

In section 3  we recall the notion of an analytic standard fan as
developed in  \cite{assi-castro-granger-jpaa-2}, and we sketch the
easy generalisation which we need from the cyclic case treated in
\cite{assi-castro-granger-jpaa-2} to the general case.

In section 4 we finally prove that the analytic standard fan is an
adapted fan in the sense of \cite{Sabbah}.

\vspace{.6cm}

\noindent Acknowledgements.- During the preparation of this paper
both authors have been partially supported by Acci\'{o}n
Integrada-Picasso HF2004-0117. First author has been  also
partially supported by MTM2004-01165. Both authors are grateful to
the D\'epartement de Math\'ematiques et Applications of the
\'Ecole Normale Sup\'erieure (Paris) where the first was invited
during the preparation of the final part of this work. They also
thank Jenifer Granger for her proofreading.

\section{Multifiltrations, Rees  rings and Rees  modules}\label{multif_rees}

Let us denote by $\cO=\CC\{x_1,\ldots,x_n\}$ the complex
convergent power series  ring and by $\cD$ the ring of germs of
linear differential operators with holomorphic coefficients (i.e.
$\cD = \cO[\partial]= \cO[\partial_1,\ldots,\partial_n]$ where
$\partial_i$ is the partial derivative with respect to $x_i$ and
the product in $\cD$ is defined by the Leibnitz's rule:
$\partial_i a=a\partial_i +
\partial_i(a)$, for each $a\in \cO$).

An operator $P$ in $\cD$ can be written as
$$P=\sum_{\beta\in \NN^n}
p_\beta(x)\partial^\beta = \sum_{\alpha,\beta\in \NN^n}
p_{\alpha\beta}x^\alpha\partial^\beta$$ where the first sum is
finite, $\alpha = (\alpha_1,\ldots,\alpha_n),
\beta=(\beta_1,\ldots,\beta_n)$, $x^\alpha=x_1^{\alpha_1}\cdots
x_n^{\alpha_n}$, $\partial^\beta=\partial_1^{\beta_1}\cdots
\partial_n^{\beta_n}$, $p_\beta(x) \in \cO$ and $p_{\alpha\beta}\in \CC$.

For each $i=1,\ldots, n$, let us remember that the $V$-filtration
on $\cD$, with respect to the hypersurface $x_i=0$, was
defined by Malgrange and Kashiwara as:
$$ V^{(i)}_\ell(\cD)= V^{(i)}_\ell = \{P=\sum_{\alpha,\beta\in \NN^n}
p_{\alpha\beta}x^\alpha\partial^\beta \in \cD \, \vert \,
\beta_i-\alpha_i \leq \ell \} $$ for each $\ell\in\ZZ$. The family
$(V^{(i)}_\ell)_{\ell \in \ZZ}$ is an increasing exhaustive
filtration on the ring $\cD$. For $i=1$, the associated graded
ring $$\gr^{V^{(1)}}(\cD) = \bigoplus_\ell
\frac{V^{(1)}_\ell}{V^{(1)}_{\ell -1}}$$ is isomorphic to the ring
$\CC\{x_2,\ldots,x_n\}[x_1,\partial_1,\partial_2,\ldots,\partial_n]$,
graded by the so-called $V^{(1)}$--graduation, where the
homogeneous elements of $V^{(1)}$--degree $\ell$ are
$$\sum_{\alpha, \beta \in \NN^n; \, \beta_1-\alpha_1=\ell} p_{\alpha
\beta} x^\alpha
\partial^\beta .$$ We have similar descriptions for $i=2,\ldots,k$.

We can also consider a $V^{(i)}$-filtration on the free module
$\cD^r$ just by defining $V^{(i)}_\ell (\cD^r)= (V_\ell^{(i)}
\cD)^r$ and more generally, for any vector
$\un{m}=(m_1,\ldots,m_r)\in \ZZ^r$, one can define the shifted
filtration $V^{(i)}[\un{m}]_\ell (\cD^r) = \oplus_{j=1}^{r}
V^{(i)}_{\ell-m_j}\cD$. Such a filtered module is called a
$V^{(i)}$-filtered free module. All the $\cD$-modules considered
will be left modules.

\begin{definicion}\label{goodfiltration}
Let $\cM$ be a $\cD$-module.  We say that a filtration $U^{(i)}(\cM)$ indexed by ${\bf Z}$ is
a good $V^{(i)}$-filtration if there exists a presentation $\cM
=\frac{\cD ^r}{\cN}$ as a quotient by a left submodule $\cN$ of $\cD ^r$, and a weight vector
$\un{m}$ such that $U^{(i)}_\ell (\cM)=\pi (V^{(i)}[\un{m}]_\ell (\cD^r)) $ where $\pi $ is
the projection $\cD ^r \to \cM$.
\end{definicion}

\subsection{$V$-multifiltrations and $V^\Gamma$-multifiltrations}\label{Vfiltrations}
Let us fix an integer  $k$ such that $1\leq k \leq n$.

For each $s=(s_1,\ldots,s_k)\in \ZZ^k$, we shall denote
$V_s(\cD)= \bigcap_{i=1}^k V^{(i)}_{s_i}(\cD)$. The family
$\{V_s(\cD)\}_{s\in \ZZ^k}$ defines a multi-filtration on the ring
$\cD$. To simplify we shall say that $V_\bullet(\cD)$ is a
$k$--filtration on  $\cD$ or even, if no confusion is possible, a
filtration on $\cD$.

Let us consider a rational simplicial cone $\Gamma$ in the first
quadrant of the dual space $(\QQ^k)^*= Hom_\QQ(\QQ^{k},\QQ)$. We
denote by $\check{\Gamma}$ the dual cone of $\Gamma$, i.e.
$$\check{\Gamma} = \{a\in \QQ^{k}\,\vert \, \gamma(a) \geq 0,
\forall \gamma \in \Gamma\}.$$
We associate
with such a cone $\Gamma$ the affine variety, denoted by $S_\Gamma$,
with coordinate ring equal to $\CC[\check{\Gamma}\cap \ZZ^{k}]$.
We will denote $A_\Gamma =\CC[\check{\Gamma}\cap \ZZ^k]$ and
$A=\CC[\NN^k]$. We denote by $\cL(\Gamma)$ the set of primitive
elements in the 1-skeleton of $\Gamma$.

The multifiltration $V^\Gamma$ on $\cD$ is defined as follows: For
each $s\in \ZZ^k$ we define $$V^\Gamma_s(\cD) = \sum_{\sigma\in
\ZZ^{k}\,\vert\, L(\sigma)\leq L(s); \forall L\in \cL(\Gamma)}
V_\sigma(\cD).$$ Notice that the sum is indexed by $\sigma \in
s-\check{\Gamma}$ and that we have the inclusion $V_s(\cD) \subset
V_s^\Gamma(\cD)$. The family $V^\Gamma_\bullet(\cD)$ is a
multifiltration of the ring $\cD$, indexed by $s\in \ZZ ^k$. This
means that~:
$$V^\Gamma_{s}(\cD).V^\Gamma_{s'}(\cD)\subset V^\Gamma_{s+s'}(\cD)
\mbox{ and } \bigcup V^\Gamma_s(\cD)=\cD.$$

We may define, in a similar way to the case of one filtration, the notion of a free multi-filtered
module and that of a good multifiltration of a  finitely generated $\cD$-module   $\cM$.
For that purpose we chose   a shift
multivector $\un {n}=(n^{(1)},\cdots  ,n^{(r)})\in ( \ZZ ^k)^r$ called also a shift matrix, with columns
$n^{(i)}\in \ZZ ^k$, and a presentation $\cM =\frac{\cD ^r}{\cN}$ of $\cM$.

\begin{definicion}
The multifiltered free module associated with $\un{n}$ is the module  $\cD ^r$
endowed with the multifiltration indexed by  $s\in \ZZ^k$, and defined as:
$$V[\un {n}] _s(\cD^r)= \bigoplus _{i=1}^r  V_{s- n ^{(i)}}(\cD  )$$
\end{definicion}

\begin{definicion} A good multifiltration  of $\cM$, is a filtration indexed by $s\in \ZZ^k$, of the type
$$U _s(\cM)= \pi (V[\un {n}]_s(\cD ^r))=\frac{V[\un{n}]_s(\cD ^r) +
\cN }{\cN}$$
for some presentation $\pi :\cD ^r \to \cM$ of $\cM.$
\end{definicion}

In the obvious sense these multifiltrations  are compatible with
the multifiltered structure on the ring $\cD$. Remark that the
$i$-th generator $e_i$ of  $\cD ^r$ is then of multidegree $n
^{(i)}\in \ZZ^k$. We may also endow $\cN $ with the induced
filtration $U_s(\cN)=V[\un{n}]_s(\cD ^r)\cap \cN$, so that we also
have $$U _s(\cM)=\frac{V[\un{n}]_s(\cD ^r) }{U_s(\cN)}.$$

We may observe that the multifiltration on the free module $\cD ^r$ is
defined as the intersection of the $V^{(i)}$--filtration with
respect to the row vectors of $\un{n}$,  $\un{n} _i=(n^{(1)}_i
,\cdots , n^{(r)}_i)\in {\bf Z}^r$, that is
$$V[\un {n}] _s(\cD^r)=V ^{(1)}[\un{n} _1]_{s_1}(\cD ^r)\cap \cdots
 \cap V ^{(k)}[\un{n} _k]_{s_k}(\cD ^r),$$
but the analogue for $\cM$, with respect to the good $V^{(i)}$-filtrations as defined
in definition \ref{goodfiltration} is not true since the inclusion:
$$U_s(\cM) = \frac{V ^{(1)}[\un{n} _1]_{s_1}(\cD ^r)\cap \cdots
 \cap V ^{(k)}[\un{n} _k]_{s_k}(\cD ^r) + \cN}{ \cN}  \subset
\bigcap _{i=1}^k \frac{ V ^{(i)}[\un{n} _i]_{s_i}(\cD ^r) +  \cN}{ \cN}
=\bigcap _{i=1}^k U^{(i)}_{s_i}(\cM)$$
may be strict.

For each good filtration there is an associated $\Gamma
$-filtration compatible with the multifiltration
$V^\Gamma_\bullet(\cD)$ on $\cD $~: \\  For each $s \in \ZZ^k$ let
us consider $$U^\Gamma_s(\cM )= \sum_{\sigma\in \ZZ^{k}\,\vert\,
L(\sigma)\leq L(s); \forall L\in \cL(\Gamma)} U_\sigma (\cM).$$
The multifiltration $U^\Gamma_\bullet(\cM)$ is a good filtration
with respect to $V^\Gamma_\bullet(\cD)$.

The goodness means here that we can verify, with the presentation
of $\cM $ as above, that $U^\Gamma_s(\cM)=\pi (V[\un {n}]^\Gamma
_s(\cD^r)) $, is still the quotient of a filtration on the free
module $\cD^r$ which is a direct sum of convenient shifts of the
filtration $V_\bullet^\Gamma $ on the ring $\cD$. More precisely
the involved filtration $V[\un {n}]_\bullet^\Gamma $ is defined by
$$V[\un {n}]^\Gamma _s(\cD^r)=
\sum_{\sigma\in \ZZ^{k}\,\vert\, L(\sigma)\leq L(s); \forall L\in
\cL(\Gamma)}V[\un {n}] _{\sigma }(\cD^r) (\,=  \bigoplus _{i=1}^r
( V^\Gamma _{s- n ^{(i)}}(\cD  ))$$ and it remains to be remarked that:
$$U_s^{\Gamma }(\cM) =\sum_{\sigma\in \ZZ^{k}\,\vert\, L(\sigma)\leq L(s); \forall L\in
\cL(\Gamma)}U _{\sigma }(\cM)=\sum_{\sigma\in \ZZ^{k}\,\vert\,
L(\sigma)\leq L(s); \forall L\in \cL(\Gamma)}\frac{V[\un
{n}]_{\sigma }(\cD ^r) + \cN }{\cN} =\frac{V[\un {n}]^\Gamma
_s(\cD^r) + \cN }{\cN}$$ so that the goodness of the filtration $U
_{\bullet }(\cM )$ implies the goodness of $U _\bullet ^{\Gamma}
(\cM )$.

\subsection{Rees  rings and Rees  modules. }\label{Reesrings}

\subsubsection {Definition of Rees rings $\cR_V(\cD)$, $\cR^\Gamma(\cD)$ and of related  Rees
modules}\label{R_VD} The Rees  ring associated with the
$V$-multifiltration on $\cD$ is defined by: $$ \cR_V(\cD) =
\bigoplus_{s\in \ZZ^k} V_s(\cD)u^s$$ where $u=(u_1,\ldots,u_k)$
are variables and the product in the Rees  ring is induced by the
natural product in the Laurent polynomial ring $\cD[u^{\pm 1}] =
\cD[u_1,u_1^{-1},\ldots,u_k,u_k^{-1}]$.

By definition the Rees  ring $\cR_V(\cD)$ is a graded
$\CC$-algebra with values group $\ZZ^k$, whose homogeneous
elements are $Pu^s$ for $P\in V_s(\cD)$ and $s\in \ZZ^k$.

Similarly given a $\cD$-module $\cM= \cD
^r/\cN$ and a shift matrix $\un{n}$, we define from the good filtrations
$V[\un{n}]_s(\cD ^r)$, and $U_s(\cM) $, the Rees modules
$$\cR_{V[\un{n}]}(\cD^r) =
\oplus_{s\in \ZZ^k} V[\un{n}]_s(\cD^r)u^s$$
$$\cR_U(\cM)= \bigoplus_{s\in \ZZ^k} U_s(\cM)u^s.$$
Both have a  natural structure of mutigraded left $\cR_V(\cD)$-module
and we shall consider $\cR_U(\cM)$ as a sub-group of $\cM[u^{\pm
1}]$, the Laurent polynomials with coefficients in $\cM$.


As $\cN $ is a left submodule of  $\cD ^r$, we can also consider
on $\cN $ the induced $V$-filtration:
$$V_s(\cN ):=\cN \cap V_s[\un{n}](\cD ^r)$$ for each $s\in \ZZ^k$ as defined before (see \ref{Vfiltrations}).
The abelian group $$\cR_V(\cN ):=\bigoplus_{s\in \ZZ^k} V_s(\cN
)u^s$$ is in fact an homogeneous submodule of $\cR_{V[\un{n}]}(\cD
^r)$. The $\cR_{V}(\cD)$--module $\cR_U(\cM)$ is naturally
isomorphic to the quotient $\cR_{V[\un{n}]}(\cD ^r)/\cR_V(\cN)$.
\begin{definicion}
We call the module $\cR_U(\cM)$ the Rees module associated with
the filtration $U_\bullet(\cM)$.
\end{definicion}

Let us now consider the filtration
$V^\Gamma_\bullet(\cD)$ as in \ref{Vfiltrations}. Then the Rees
ring associated with this filtration is defined in a similar way as
$$\cR^\Gamma(\cD) =
\bigoplus_{s\in \ZZ^ k} V_s^\Gamma(\cD) u^s \subset \cD[u^{\pm
1}]$$ the product being induced by the one of $\cD[u^{\pm 1}]$.

If $\cM=\frac{\cD ^r}{\cN}$ is a finitely
generated $\cD$-module and $U_\bullet (\cM)$ is a good filtration
on $\cM$ with respect to $V_\bullet(\cD)$ then, for any cone
$\Gamma$ as in \ref{Vfiltrations}, the abelian group
$$\cR^\Gamma(\cM) = \bigoplus_{s\in \ZZ^k} U^\Gamma_s(\cM) u^s
\subset \cM[u^{\pm 1}]$$ is a left graded
$\cR^{\Gamma}(\cD)$-module.
\begin{definicion}
We call the module $\cR^\Gamma(\cM)$
the Rees module associated with the filtration $U_\bullet(\cM)$
and the cone $\Gamma$.
\end{definicion}

The fact that a filtration is good in the sense of previous
section \ref{Vfiltrations} is then equivalent to the fact that the
Rees module of $\cM$ is finitely generated over $\cR_V(\cD)$.
Similarly the goodness of $U^\Gamma_\bullet(\cM)$ is equivalent to
the fact that the Rees module $\cR ^{\Gamma}{(\cM)}$ as defined
above, is finitely generated over $\cR^{\Gamma}{(\cD)}$. Since
$\cR^{\Gamma}{\cD}$ is a Noetherian ring, this implies in
particular an Artin type property: If $\cM' \subset \cM$ is a
submodule of a $\cD$-module $\cM$ endowed with a good filtration
$U_\bullet(\cM)$ then the induced filtration on $\cM'$, namely
$U^\Gamma_\bullet(\cM)\cap \cM'$ is good.

The ring $\cR^{\Gamma}{(\cD)}$ contains as a subring the
$\CC$--algebra $A_\Gamma = \CC[\check{\Gamma}\cap \ZZ^k]\subset
\CC [u,u^{-1}]$, all these rings being included in
$\cD[u,u^{-1}]$.

\begin{teorema}\label{forme_faible} There is a fan $\cE$ in $(\QQ^k)^*_+$ such that for each
cone $\Gamma\in\cE$ the Rees module $\cR^\Gamma(\cM)$ is flat over
$A_\Gamma$.
\end{teorema}

We will prove a more precise form of this  result we in
section \ref{flatness}, see theorem \ref{mainresult}, after having
recalled in section \ref{analytic_standard_fan} the notion of an
analytic standard fan.

\subsubsection{The fiber at $0$}\label{fibres} Let us denote by $\mathfrak m$ the ideal of $\cR_{V}(\cD)$
generated by $(u_1,\ldots,u_k)$.  This is a two sided ideal since
it is generated by central elements. The fiber at zero of
$\cR_{V}(\cD)$ (resp. of $\cR_{V[\un{n}]}(\cD ^r)$)  is by
definition the quotient ring (resp. the quotient module),
$$\frac{\cR_{V}(\cD)}{\mathfrak m} \,\,(\mbox{resp. }
\frac{\cR_{V}[\un{n}](\cD ^r)}{\mathfrak m.\cR_{V}[\un{n}](\cD
^r)}).$$ More generally the fiber at zero of
$\cR_{U}(\cM)=\frac{\cR_{V}[\un{n}](\cD ^r)}{\cR_{V}(\cN)}$ is the
quotient $$\frac{\cR_U(\cM)}{{\mathfrak m }{\cR_U(\cM)}}$$ which
is naturally isomorphic to the quotient
$$\frac{\cR_{V[\un{n}]}(\cD^r)}{\cR_{V}(\cN )+\mathfrak m.\cR_{V[\un{n}]}(\cD^r )}.$$

Notice that the fiber at 0 of the module $\cR_U(\cM)$ can  be zero
for a non-zero $\cR_U(\cM)$, as shown in the following example.
Let us denote by $I$ the principal ideal of $\cD$ generated by the
differential operator $P=1+x_1^2\partial_1$ and let us suppose
$k\geq 1$. Then $\cR_{V}(I)$ is the principal ideal of
$\cR_{V}(\cD)$ generated by $Pu_1^0=1+(x_1^2\partial
_{1}u_1^{-1})u_1$ and the fiber at zero of $\cR_{V}(\cD/I)$ is
then zero since $\cR_{V}(I )+\mathfrak m.\cR_{V}(\cD)$ contains
$P-(x_1^2\partial_{1}u_1^{-1})u_1=1$, so that it is equal to
$\cR_{V}(\cD).$

Similarly for any $k$-dimensional cone $\Gamma$ the $\CC$--algebra
$A_\Gamma = \CC[\check{\Gamma}\cap \ZZ^k]$ has a maximal ideal
$$\mathfrak m^{\Gamma }= \CC[\check{\Gamma}\cap \ZZ^k\setminus\{ 0\}]$$
because $\check{\Gamma}$ is strictly convex. Then we define the
fiber at zero of the Rees module $\cR^\Gamma(\cM)$ as
$\frac{\cR^\Gamma(\cM)}{\mathfrak m^{\Gamma }\cR^\Gamma(\cM)}.$ As
is explained in \cite{Sabbah}, this fiber is a module over the
ring $\frac{\cR^\Gamma(\cD)}{\mathfrak m^{\Gamma
}\cR^\Gamma(\cD)}\simeq gr ^{\Gamma}\cD$.

\vspace{1cm}

\subsubsection{Description of  Rees rings $\cR_V(\cD)$ and $\cR^\Gamma(\cD)$}\label{RgammaD}

It is useful to describe the Rees  ring $\cR_{V}=\cR_V(\cD)$
(resp. $\cR^\Gamma(\cD)$) as subrings   of the  ring of relative
differential operators $\cD_{\CC^n\times \CC^k/\CC^k}$ (resp
$\cD_{\CC^n\times S_{\Gamma}/S_{\Gamma}}$). This is nothing but an
explicit version of the interpretation of a Rees ring as a ring of
relative differential operators on the deformation of $Y\times
0_{\Gamma}$ to its normal cone in $\CC^n\times S_{\Gamma},$ see
\cite {Sabbah}. We denote by $S_\Gamma$ the algebraic variety
associated to the ring $A_\Gamma$.

Let us define first
$$\cA = \CC\{X'U,X''\}[X',U,\Delta]=
\CC\{X_1U_1,\ldots,X_kU_k,X_{k+1},\ldots,X_n\}[X',U,\Delta]$$
where $U=(U_1,\ldots,U_k)$, $X'=(X_1,\ldots,X_k)$,
$X''=(X_{k+1},\ldots,X_{n})$, $X=(X_1,\ldots,X_n)$,
$\Delta=(\Delta_1,\ldots,\Delta_n)$ are new variables satisfying
the following relations, for $i=1,\ldots,n$:
$$\Delta_i X_i = X_i \Delta_i + 1$$ the other relations being
trivial.

The ring $\cA$ is graded with $\ZZ^k$ as a values group (we will say
that $\cA$ is a $\ZZ^k$--graded ring) , namely we have $$\cA =
\bigoplus_{s\in \ZZ^k} \cA_s$$ where $\cA_s$ is the set $$ \cA_s =
\{\sum_{\alpha,\beta\in \NN^n; \sigma \in \NN^k} f_{\alpha\beta
\sigma} X^\alpha \Delta^\beta U^\sigma\in \cA \,\vert \,
f_{\alpha\beta\sigma}\in \CC\{X'U,X''\}, \sigma_i+\beta_i-\alpha_i
= s_i; i=1,\ldots,k\}.$$

\begin{propo}
There is an isomorphism  of $\ZZ^k$--graded rings
$$i=i_{V} : \cR_{V}(\cD) \rightarrow \cA$$ defined by:

\begin{itemize}
 \item $i(u_j)=U_j$ for $j=1,\ldots,k$.
\item $i(x_ju_j^{-1}) = X_j$ for $j=1,\ldots,k$ \; \mbox{ and}
\quad
$i(x_j)=X_j$ for $j=k+1,\ldots,n$.
 \item $i(\partial_ju_j)=
\Delta_j$ for $j=1,\ldots,k$
 \; \;  \mbox{ and} \quad $i(\partial_j)=
{\Delta_j}$ for $j=k+1,\ldots,n$.
\end{itemize}
\end{propo}

\begin{prueba}{.}
It is clear that $i$ is injective and that, by the formula
$i(x^{\alpha}\partial ^{\beta}u^s)=X^{\alpha}\Delta ^{\beta}(\Pi
_{i=1}^kU_i^{s_i+\alpha_i-\beta_i}),$
  $i(V_s(\cD)u^s)= \cA_s$ for all $s\in \ZZ^k$ .
\end{prueba}

Let us give now the same description for the ring $\cR^\Gamma(\cD)$. We shall
do it in the only case of interest for us,
when the cone $\Gamma$ is basic, which means that if it is defined by $k$
independent linear forms $\{L_1,\ldots,L_k\}$ generating the
lattice $\ZZ^k$. For the sake of simplicity
let us suppose $k=n$.

In this situation we are going to make a change
of variables in order to write $A_\Gamma $  as a polynomial ring.

For any $i$ we can write $L_i=(\ell_{i1},\ldots,\ell_{ik})$ and we
can suppose by a suitable ordering of the linear forms $L_i$  that
the determinant of the matrix $\cL=(\ell_{ij})$ equals  1.

Let us write $\cL'$ for the inverse matrix of $\cL$, then the
columns $\{C _1,\ldots,C _k\}$ of $\cL'$ form a basis of the dual
cone $\check{\Gamma}$.

 The family $u^{C _1},\ldots,u^{C _k}$ generates the
$\CC$--algebra $A_\Gamma = \CC[\check{\Gamma}\cap \ZZ^k]$. Here we
write $u^{C_j}=u_1^{c_{1j}}\cdots u_k^{c_{kj}}$ where the
$c_{ij}$'s are the entries of the column $C_j$.

We will write $W_i=U^{C_i}$ for $i=1,\ldots,k$. Similarly, we have
$U_j = W^{C_j(\cL)}$ for $j=1,\ldots,k$ where $C_j(\cL)$ is the
$j-th$ column of the matrix $\cL$. In section \ref{flatness} we
will only use the column $C_j=C_j(\cL')$. \\
Let us remark here
that if $k<n$ then we can still define $W$, the relationship
between $W$ and $U$ being exactly the same.

\begin{propo}
The ring $\cR^\Gamma(\cD)$ is isomorphic to the subring of
$\cA'=\CC\{X,W\}[\Delta]$, which contains all polynomials and in
which only convergent power series with respect to the
family of monomials $X_iU_i=X_iW^{C_i(\cL)}$ for $i=1,\cdots k$,
and $X_{k+1},\cdots , X_n$, in the general case $k\leq n$  are allowed.
\end{propo}
\begin{prueba}{.}
A monomial $X^\alpha \Delta^\beta U^{s-\beta+\alpha}$ can be
written in terms of $X,\Delta,W$ as $X^\alpha \Delta^\beta
W^{\cL(s-\beta+\alpha)}.$ All that  remains to be done is to enumerate
the monomials derived from
 the power series variables $x_1,\cdots , x_n$
\end{prueba}

 The ring $A_{\Gamma}$ is identified to the subring
$\CC[W](=\CC[W_1,\ldots,W_k]) \subset \cR^\Gamma(\cD)$ and this inclusion is flat. We
shall prove this statement in \ref{RGammaDplat}.

The fiber of $\cR^\Gamma(\cD)$ at the origin is by definition
$\cR^\Gamma(\cD)\otimes_{\CC[W]} \frac{\CC[W]}{(W)}$ which is
isomorphic to a Weyl algebra, namely the Weyl algebra
$\CC[X,\Delta]$. This Weyl algebra is endowed with a
$\ZZ^{k}$--graduation by $weight(X_i)=-\epsilon_i$ and
$weight(\Delta_i)=\epsilon_i$ for $i=1,\ldots,k$ where
$\epsilon_i$ is the vector in $\ZZ^{k}$ whose $j-th$ coordinate is
$\delta_{ij}$.

\subsection{Multifiltrations and Bernstein-Sato functional
equations}\label{bfunctions}

Let us detail the functional equation problem raised in the
introduction. It has already been remarked that it gives rise to a
situation where a multifiltration along transverse hypersurfaces
comes out in a natural way.  These equations are of the type
$$P(\lambda)f^{\lambda }=b(\lambda )f_1^{\lambda  _1+1}\cdots
f_k^{\lambda  _k+1}$$ where $\lambda =(\lambda _1,\cdots , \lambda
_k)$ is a $k$-uple of indeterminates, and $f=(f_1,\cdots , f_k)$
is a $k$-uple of analytic functions on $X$. They  are naturally
written in the module $\cO _X\big[ \lambda _1,\cdots , \lambda
_k,\frac{1}{f_1 \cdots f_k} \big]f^{\lambda }$, endowed with a
structure of a $\cD _X[\lambda ]$-module, which can be extended in
a way discovered by B. Malgrange to a $\cD _{X\times \CC
^k}$-module structure with an action of the $2k$ variables $t_j,
\partial _{t_j}$such that $\lambda _j=-\partial _{t_j}t_j.$ The module
$\cM$ that we have then to consider is the module generated over
the ring $\cD _{X\times \CC ^k}$ by $f^{\lambda }$, with its
naturally defined multifiltration $V_{\bullet}(\cD _{X\times \CC
^k})\cdot f^{\lambda }$ along the hypersurfaces $t_1=0, \cdots ,
t_k=0$. The proof of Sabbah in \cite{Sabbah} can be sketched as
follows: let us define for any linear  form $L\in (\QQ
^k)^{\star}$ with positive rational coefficients $\ell _j \geq 0$,
a filtration of $U^L_{\bullet}(\cM)$, associated with the linear
combination $\sum \ell _jV^{(j)}$ of the basic $V$- filtrations
(see notations at \ref{multif_rees}). Precisely we set (see also
the notation \ref{filtrationVL} below):
$$
V_{\bullet}^{L}(\cD _{X\times \CC ^k})=\{ P=\sum_{\mu,\nu \in
\NN^k} P_{\mu\nu}(x,\partial _x)t^{\mu} \partial _t^{\nu} \,, \;
\sum \ell _i(\nu _i-\mu _i)\leq \bullet \}$$
$$ U^L_{\bullet}(\cM)=V_{\bullet}^{L}(\cD _{X\times \CC ^k})\cdot
f^{\lambda }$$
 Then the existence of a functional equation, in
which $b(\lambda )$ is a product of affine forms $L(\lambda )+c$
comes out from the following three steps:

1) There are Bernstein-Sato polynomials $b_L(\lambda )$, relative
to each $L$, with functional  equations
$$b_L(\lambda )f^{\lambda }\in V_{<0}^L(\cD _{X\times \CC
^k})\cdot f^{\lambda }.$$

2) The saturated filtration
$$\ov{U_s}(\cM)=\bigcap _{L\in (\QQ ^k)^{\star}}U_{L(s)}^L(\cM)$$
can be defined by using only a {\em finite number} of fixed linear
forms.

3) The saturated filtration is good which is equivalent to the
existence  of a $k$-uple of integers $\kappa $, such that
$$\forall s\in \NN ^k,  \;\; U_s(\cM)\subset \ov{U_s}(\cM)\subset U_{s+\kappa }(\cM).$$

It is step 2)  which, in the proof in \cite{Sabbah}, makes an essential
use of  the notion of an adapted
fan, whose existence is proved in this paper, see theorem
\ref{forme_faible}. In  Bahloul's paper
\cite{bahloul}, this step is made by constructive methods which do not use the
flatness property.

\section{The analytic standard fan}\label{analytic_standard_fan}

In this section we will summarize the main results of
\cite{assi-castro-granger-jpaa-2}. Since these results are only
given for a module $\cD/I$  over $\cD$ which is a quotient by an
ideal we will then show briefly how to adapt them to a module of
the type $\frac{\cD ^r}{\cN}$. Let us remark that the fan used in
this paper is obtained by a restriction of the fan in
\cite{assi-castro-granger-jpaa-2} to a subset of linear forms of
the set $\cU$ defined below, namely the linear combinations of the
filtrations $V^{(i)}$.

Let $\cU$ be the set of linear forms $\Lambda: \RR^{2n}
\rightarrow \RR$, $\Lambda(\alpha,\beta)=\sum_{i=1}^ n e_i
\alpha_i + \sum_{i=1}^ n f_i \beta_i$ with $e_i+f_i \geq 0$ and
$e_i \leq 0$ for $i=1,\ldots,n$. If $$P=\sum_{\alpha \beta}
p_{\alpha \beta} x^\alpha
\partial^\beta$$ is an element in $\cD$ we define $\ord^\Lambda(P)$ --the
$\Lambda$-order of $P$-- to be the maximal value of
$\Lambda(\alpha,\beta)$ for $\alpha,\beta$ such that $p_{\alpha
\beta} \not= 0$.

The $\Lambda$-filtration $F_{\Lambda,\bullet} (\cD)$ is defined by
$$F_{\Lambda,\ell} = F_{\Lambda,\ell}(\cD) = \{P \in \cD \,\vert \, \ord^\Lambda(P)\leq
\ell\}$$ for any $\ell \in \RR$. We will write
$F_{\Lambda,<\ell}(\cD):= \{P \in \cD \,\vert \, \ord^\Lambda(P) <
\ell\}$. If $e_i=0, f_i=1$ for $i=1,\ldots,n$ the corresponding
$\Lambda$--filtration is nothing but the usual filtration by the
order of differential operators. We shall denote it simply
$F_\bullet(\cD)$.

Let us recall that we have fixed $k\leq n$. For each linear form
$L\in (\QQ^k)^*_+$ (i.e. the coefficients of $L$ are non-negative)
we denote by $\widetilde{L}$ the linear form on $\RR^{2n}$ defined
by $\widetilde{L}(\alpha,\beta)=L(\beta)-L(\alpha)$ where
$L(\alpha_1,\ldots,\alpha_n)=L(\alpha_1,\ldots,\alpha_k)$.

\begin{notacion}\label{filtrationVL}
We define the filtration $V_\bullet ^L$ on $\cD$, indexed by the
set of values $L(\ZZ^k)$, as:
$$V^L_{\ell}(\cD)=F_{\widetilde{L},\ell}(\cD)=\{ P\in \cD \,\vert \, \ord^{\widetilde{L}}(P)
\leq  \ell\}$$
\end{notacion}


Let us fix $L\in (\QQ^k)^*_+$.  The graded ring  associated to the
filtration $V^L=F_{\widetilde{L}}$ on $\cD$,
is by definition
$$\gr^{\widetilde{L}}(\cD) = \bigoplus_{\ell \in \widetilde{L}(\ZZ^{2n})}
\frac{F_{\widetilde{L},\ell}}{F_{\widetilde{L},<\ell}}.$$ If no
confusion is possible we shall write simply $\ord^L$ and
$\gr^L(\cD)$ instead of $\ord^{\widetilde{L}}$ and
$\gr^{\widetilde{L}}(\cD)$.

The graded ring $\gr^L(\cD)$ is a ring of differential operators
and its structure is the following: suppose the coefficients of
the form $L$ are $(e_1,\ldots, e_k)$ and suppose we also have
ordered the variables to have $e_i>0$ for $1 \leq i \leq \ell$ for
some $\ell \leq k$.

Then the graded ring $\gr^L(\cD)$ is isomorphic to the ring
$$\CC\{x_{\ell+1},\ldots,x_n\}[x_1,\ldots, x_{\ell},
\partial_1,\ldots,\partial_{n}].
$$
In this ring the graduation is induced by the weights
$\weight(x_i)=-e_i, \,\weight(\partial_i)=e_i$ for $1\leq i \leq
\ell$ and $\weight(x_j)=\weight(\partial_j)=0$ otherwise. There is
only a finite number of types of these  rings, one for each
partition of $\{1,\ldots,k\}$ into two sets.

\begin{notacion} \label{sigmaL} For each $P\in \cD$ and for each $d \in L(\ZZ^{k})$ with
$\ord^L(P) \leq d$ we denote by $\sigma^L_d(P)$ --the $L$--symbol
of $P$ of order $d$-- the class of $P$ in $F_{\wL,d}/F_{\wL,<d}$.
The principal symbol of $P$ is by definition
$\sigma^L(P)=\sigma^L_{d}(P)$ if $d=\ord^L(P)$. For $P,Q \in \cD$
we have $\sigma^L(PQ) = \sigma^L(P)\sigma^L(Q)$. For any left
ideal $I$ in $\cD$ we denote by $\gr^L(I)$ the graded ideal of
$\gr^L(\cD)$ generated by the set $\{\sigma^L(P)\,\vert\, P \in
I\}$. We set $\sigma^L(0)=0$.
\end{notacion}

We denote by $\cD[t]$ the $\CC$-algebra $\cO[\partial,t] =
\CC\{x_1,\ldots,x_n\}[\partial_1,\ldots,\partial_n,t]$ with
relations ($t$ being a new variable)
\begin{itemize}
\item $[t,a]=[t,\partial_i]=[a,b]=[\partial_i,\partial_j]=0$,
\item $[\partial_i,a] = \frac{\partial a}{\partial x_i}t$,
\end{itemize}
for $a,b\in \cO$ and $i=1,\ldots,n$.

The ring $\cD[t]$ is isomorphic to the Rees ring associated with
the order filtration $F_\bullet$ on $\cD$. Since this Rees ring is
by definition $$\cR_F(\cD) = \bigoplus_{\ell\in \ZZ}
F_\ell(\cD)v^\ell \subset \cD[v]=\cD\otimes_\CC \CC[v]$$ for a new
variable $v$, we can define an isomorphism of graded rings $\iota
: \cD[t] \rightarrow \cR_F(\cD)$ by $\iota(a)=a=av^0,
\iota(\partial_j) =
\partial_j v,$ and $\iota(t)=1v=v$.

The natural graded structure of $\cR_F(\cD)$ can be translated on
$\cD[t]$. An homogeneous element of degree $d\in \ZZ$ in $\cD[t]$
is nothing but an expression $$\sum_{\ell + \vert \beta \vert = d}
a_{\ell\,\beta}\partial^\beta t^\ell$$ for some
$a_{\ell\,\beta}\in \cO.$

For $P=\sum_\beta p_\beta(x) \partial ^\beta \in \cD$, the element
$Pv^{\ord(P)}\in \cR_F(\cD)$ is called the homogenization of $P$,
where $\ord(P)$ is the usual order of $P$. It is useful to see
$Pv^{\ord(P)}$ as an element of $\cD[t]$. The homogenization of
$P$ is then denoted by $h(P)$ and we have
$$h(P)= \sum_\beta p_\beta(x)\partial^\beta t^{d-|\beta|}$$ for
$d=\ord(P)$ the usual order of $P$.

For each linear form $L\in (\QQ^k)^*_+$ we can define in a natural
way a filtration $V^L_\bullet(\cD[t])$ on $\cD[t]$, the $L$--order
--denoted $\ord^L(R)$-- of an element $$R=\sum_{\ell\alpha \beta }
r_{\ell \alpha \beta } x^\alpha\partial^\beta t^\ell$$ being the
maximal value of $L(\beta)-L(\alpha)$ for $r_{\ell\alpha\beta}
\not= 0$. The associated graded ring is denoted by
$\gr^L(\cD[t])$. We have a natural ring isomorphism  from
$\gr^L(\cD[t])$ onto $\gr^L(\cD)[t]$, where with the notations as
before we have $[\partial_i,x_j]=\delta_{ij}t$, $\delta_{ij}$
being the Kronecker symbol.

For each left ideal $I$ in $\cD$ we denote by $h(I)$ the left
homogeneous  ideal of $\cD[t]$ generated by $\{h(P)\,\vert \, P\in
I\}$. As in the case of $\cD$, we denote by $\gr^L(h(I))$ the
homogeneous  ideal --in $\gr^L(\cD[t])$-- generated by the set of
principal symbols $\sigma^L(G)$ of the elements $G$ in $h(I)$.

The main result of \cite{assi-castro-granger-jpaa-2} is the
following

\begin{teorema} Let $I$ be a non-zero left ideal of $\cD$ and let $h(I)$ be the associated
homogenized ideal in $\cD[t]$. Then there exists a partition $\cE$
of $\cU$ into convex rational polyhedral cones such that for any
$\Gamma \in \cE$ the ideals $\gr^\Lambda(h(I))$ and
$\gr^\Lambda(I)$ do not depend on $\Lambda\in \Gamma$.
\end{teorema}

Let us sketch the generalization of this result to the case  of a
$\cD $-module $\cM =\frac{\cD ^r}{\cN}$ given by a general
presentation, and endowed with a good multifiltration
$$U_s(\cM ) =
\frac{V_s[\un{n}](\cD ^r)}{V_s(\cN)}=\frac{V_s[\un{n}](\cD ^r) +
\cN }{\cN}$$ associated with a multivector shift $\un{n}\in
(\ZZ^k)^r$.

Here we have to restrict the set  $\cU$ to the set of linear forms
$\wL$, indexed by $L$ in $(\QQ^k)_+^*$, as defined above:
$$\wL(\alpha , \beta )= L(\beta)-L(\alpha).$$

Let us recall  the notation  \ref{filtrationVL} $V_{L(s)}^{L}(\cD)
= F_{\wL,L(s)}(\cD)$, which we shall extend to any free
multifiltered module.

Precisely, given a shift multivector
$\un{n}=(n^{(1)},\ldots,n^{(r)})\in (\ZZ^k)^r$ we define a family, indexed by $L\in
(\QQ^k)_+^*$, of $L({\bf Z}^k)$-filtrations on the module $\cM$  as follows:

$$V_{\ell}^{L}[L(\un{n})](\cD ^r) =
\bigoplus _{i=1}^r V_{\ell-\ell_i}^{L}(\cD)$$ where
$\ell_i=L(n^{(i)})$ for $i=1,\ldots,r$ and
$L(\un{n})=(\ell_1,\ldots,\ell_r)$.

$$U_{\ell}^{L}(\cM) =\frac{V_{\ell}^{L}[L(\un{n})](\cD ^r)}
{ V_{\ell}^{L}(\cN)}= \frac{V_{\ell}^{L}[L(\un{n})](\cD ^r) + \cN
}{\cN }.$$

We again have the notion of a graded associated module
$$\gr^{L}(\cM) = \bigoplus_{\ell \in L(\ZZ^{k})}
\frac{ U_{\ell}^{L}(\cM)}{ U_{<\ell}^{L}(\cM) } =\frac{\gr^{L}(\cD
^r)}{\gr^{L}(\cN)}.$$

We can endow the module $\cD[t]^r$ with a filtration of a filtered
$\cD  [t]$-module by $$V_{\ell}^{L}[L(\un{n})](\cD [t] ^r) =
\bigoplus _{i=1}^r V_{\ell-\ell_i}^{L}(\cD [t])$$ with the same
notations as before,  and the homogeneization of $(P_1, \ldots
,P_r)\in \cD ^r $ is just $$(t^{d-d_1}h(P_1), \ldots
,t^{d-d_r}h(P_r))$$ where $d_i$ is the usual order of the operator
$P_i$ and $d=\max d_i$.

We can now state the straightforward generalisation of the main
theorem in  \cite{assi-castro-granger-jpaa-2}.

\begin{teorema} Let $\cN$ be a non-zero left submodule  of $\cD ^r$
endowed with a good filtration,
and let $h(\cN)$ be the associated homogenized submodule in
$\cD[t]^r$. Then there exists a partition $\cE$ of $(\QQ^k)_+^*$
into convex rational polyhedral cones such that for any $\Gamma
\in \cE$, the submodules $\gr^{L}(h(\cN))$ and $\gr^{L}(\cN)$ do
not depend on $L\in \Gamma$. \end{teorema}

The proof of the theorem is identical to the one in
\cite{assi-castro-granger-jpaa-2}, since the main ingredient,
which is the division theorem, can be adapted to the module case.
See for example \cite{granger-oaku}, where an even more
general situation (with shifts for the variables $t$) is treated. In
particular we still have on every cone of the partition $\cE$ the
notion of a reduced Gr\"obner basis of $h(\cN )$, valid for any
linear form $L$ in the cone.

\section{Analytic standard  fan and flatness}\label{flatness}

After  Herrmann-Orbanz \cite{Herrmann-Orbanz} we have the
following criterion for flatness of a module $M$ over a
commutative ring $R,$ both being graded by a commutative group $G$:
$$M=\bigoplus_{g\in G} M_g , \quad R = \bigoplus_{g\in G} R_g $$
\begin{teorema}
The (not necessarily of finite type) graded $R$--module $M$ is
$R$-flat if and only if for any graded ideal $H\subset R$, we
have:
$$Tor_1^R(M,R/H)=0.$$
\end{teorema}

We will apply this criterion to  the Rees module $$\cR ^{\Gamma
}({\cM}) = \bigoplus _{s\in \ZZ ^k}U^{\Gamma }_s(\cM) u^s$$ graded
by the group $\ZZ^k$, seen as a module over the graded ring $\cR
^{\Gamma }({\cD})$, and therefore also as a module over the graded
subring $R=\CC[W]= \CC[W_1,\cdots, W_k]$ (see \ref{Reesrings}).

Note that for any $a\in \NN^k$, $\CC \cdot W^a$, is a homogeneous
subspace of $R$ generated by a non zero element of $\ZZ^k$-degree
$\cL' a$, since with the notation of section \ref{RgammaD} $W^a =
u^{\cL'a}$. Therefore the graded ideals of $R$
we have to take into consideration are the monomial ideals of the
polynomial ring $R$.

%
%
%

We now deal with the case of $\cM=\frac{\cD ^r}{\cN}$. If we endow
$\cN$ and $\cM$ with the induced and quotient filtration, from a
$V_{\bullet}[\un{n} ]$--filtration  on the free module $\cD ^r$,
we have the exact sequence:
$$ ({\cal S}) \qquad 0 \longrightarrow \cR ^{\Gamma }({\cN })  \longrightarrow \cR
^{\Gamma }({\cD ^r})
\longrightarrow \cR ^{\Gamma }({\cM}) \longrightarrow 0.$$

\begin{lema}\label{RGammaDplat}
The ring $\cR ^{\Gamma }({\cD})$ is flat over $\CC[W]$.
\end{lema}

\begin{prueba}{.}
By the theorem of Herrmann-Orbanz, it is sufficient to prove that
$$Tor_1^{\CC[W]}(\cR^\Gamma(\cD),\CC[W]/H)=0$$ for any homogeneous
ideal $H$. If we consider the tensor product by $\cR ^{\Gamma
}(\cD )$ of the exact sequence
$$0 \longrightarrow  H  \longrightarrow \CC[W]
\longrightarrow \CC[W]/H \longrightarrow 0,$$ we see that this
vanishing result is equivalent to the injectivity of the mapping
$$H  \otimes \cR ^{\Gamma }(\cD ) \longrightarrow \cR ^{\Gamma
}(\cD).$$ \indent Since the ideal $H$ is multihomogeneous, it is a
monomial ideal and we may consider an element \break
$\mathfrak{Q}=W^{a_1}\otimes Q_1+\ldots + W^{a_r}\otimes Q_r$ of
the kernel of this map. Let $Q_i=\sum q_{i,\alpha ,\beta ,
\ell}X^{\alpha } \Delta ^{\beta}W^{\ell }$ be the development of
$Q_i$.

The fact that $\mathfrak{Q}$ is in the kernel can be expressed by the
equation
$$(\star )\qquad \sum _{i=1}^r W^{a_i}Q_i=0 .$$

Let us consider the partition of $\bigcup (a_i+\NN ^k)\subset  \NN ^k$ given
by~:
$$\Delta _1=a_1+\NN ^k\; , \quad \Delta _i=(a_i+\NN ^k)\setminus
(\Delta _1\cup \ldots \cup \Delta _{i-1}), \mbox{ for } i=2,\cdots ,r$$
We may write $Q_i=\sum _{p=1} ^i Q_{i,p}$, with~:
$$Q_{i,p}=\sum _{\ell+a_i\in\Delta _p} q_{i,\alpha ,\beta , \ell}X^{\alpha }
\Delta ^{\beta}W^{\ell },$$ operators having disjoint $(W)$-Newton
diagrams. We set $a_i+v_{i,p}=a_p+w_{i,p}$, with
$v_{i,p},w_{i,p}\in \NN ^k$ minimal for the componentwise partial
ordering, and then we may write $Q_{i,p}=W^{v_{i,p}}R_{i,p}$ with
$R_{i,p}\in \cD$, and we may also remark that $W^{a_i}Q_{i,p}$ is
the part of the operator $W^{a_i}Q_i$, whose Newton diagram  is
contained in $\Delta _p$ so that the equation $(\star )$ splits
into:
$$(\star )_p \qquad \sum _{i=p}^r W^{a_i}Q_{i,p}= W^{a_p}\sum _{i=p}^r
W^{w_{i,p}}R_{i,p}=0 \,\mbox{ for }\, p=1, \ldots ,r .$$ This
finally implies the desired result by the following calculation~:
\begin{align*}
\mathfrak{Q}=& \sum _{i=1}^rW^{a_i}\otimes
( \sum _{p=1}^i W^{v_{i,p}}R_{i,p})=\sum _{p=1}^r \sum _{i=p}^r W^{a_i+v_{i,p}}
\otimes R_{i,p} \\ =& \sum _{p=1}^r \sum _{i=p}^r W^{a_p+w_{i,p}}
\otimes R_{i,p} = \sum _{p=1}^r W^{a_p} \otimes (\sum _{i=p}^r
W^{w_{i,p}} R_{i,p})=0
\end{align*}
\end{prueba}

By the theorem of Herrmann-Orbanz, the flatness of
$\cR^\Gamma(\cM)$ -- if $\Gamma$ is included in a cone of the
analytic fan of $h(\cN)$ -- is a consequence of the equality
$$Tor_1^{\CC[W]}(\cR^\Gamma(\cM),\CC[W]/H)=0$$ for any homogeneous, that
is monomial, ideal $H\subset \CC[W]$.

Because of the above  lemma, and of the exact sequence $({\cal
S})$, this $Tor$ module appears as the kernel of the map:
$$\CC[W]/H\otimes_{\CC[W]} \cR ^{\Gamma }({\cN}) =
\frac{\cR ^{\Gamma }({\cN})}{H\cdot \cR ^{\Gamma }({\cN }) }
\longrightarrow  \frac{\cR ^{\Gamma }({\cD}^r)}{H\cdot \cR
^{\Gamma }({\cD}^r)}$$ and we will prove that it is injective.

The injectivity means that $H\cdot \cR ^{\Gamma }({\cN})= H\cdot
\cR ^{\Gamma }({\cD}^r) \cap \cR ^{\Gamma }({\cN})$ for any
monomial ideal $H = \sum _{i=1}^\ell \CC [W]W^{a_i}.$ The
inclusion $H\cdot \cR ^{\Gamma }({\cN})\subset H\cdot \cR ^{\Gamma
}({\cD}^r) \cap \cR ^{\Gamma }({\cN})$ is obvious.

\begin{teorema} \label{mainresult} Let $\Gamma$ be a cone included in a cone
of the analytic standard fan of $h(\cN)$. Then $\cR^\Gamma(\cM)$ is
flat over $A_\Gamma=\CC[W]$. \end{teorema}

\begin{prueba}{.} We have to prove the inclusion $H\cdot \cR
^{\Gamma }({\cN}) \supset H\cdot  \cR ^{\Gamma }({\cD^r}) \cap \cR
^{\Gamma }({\cN})$ for any monomial ideal $H$ in $\CC[W]$. It is
enough to treat the case  where $H$ is generated by a set of
variables $\{W_{j_1},\ldots,W_{j_p}\}$ for some
$J=\{j_1,\ldots,j_p\}\subset \{1,\ldots,k\} $. We denote by $W_J$
the ideal generated by $\{W_{j_1},\ldots,W_{j_p}\}$ and we call it
a coordinate ideal. This is due to the two following lemmata.

\begin{lema}
Let $H_1\subset H_2 \subset \CC[W]$ be two ideals and $M$ a
$\CC[W]$--module  such that
$$Tor_1^{\CC[W]}(M,\CC[W]/H_2)=Tor_1^{\CC[W]}(M,H_2/H_1)=0$$
then $Tor_1^{\CC[W]}(M,\CC[W]/H_1)=0$.
\end{lema}
\begin{lema} Let $H$ be a monomial ideal in
$\CC[W]=\CC[W_1,\ldots,W_k]$. Then there exists a sequence of
monomial ideals $$H=H_0 \subset H_1 \subset \cdots \subset
H_r=\CC[W]$$ such that for any $i$ the quotient  $H_{i+1}/H_i$ is
isomorphic to $\CC[W]/W_{J(i)}$ for some $J(i)\subset
\{1,\ldots,k\}$.
\end{lema}

\begin{prueba}{.}  By noetherianity of $\CC[W]$ it is sufficient to
construct $H_1=H_0+\CC[W]m$ where $m$ is a monomial in $\CC[W]$
such that $m\not\in H_0$ and $(H_0:m)$ is some coordinate ideal
$W_J$. In the Artinian case just take an element in the socle and
in general use an induction on the dimension of $\CC[W]/H.$
\end{prueba}

So, in order to prove theorem \ref{mainresult} it is sufficient to
treat the case $H$ generated by $\{W_1,\ldots,W_p\}$ for any
$p\leq k$. We have (see \ref{RgammaD}) $W_i = U^{C_i}$ for
$i=1,\ldots,k$ where $C_i\in \ZZ^k$ is the $i$-th column of the
matrix $\cL'$. Since the considered submodules of $\cR ^{\Gamma
}(\cD ^r ) $ are $V_{\bullet}[\un{n} ]$--homogeneous, it is sufficient
to prove that any homogeneous element $Qu^s \in H\cdot
\cR^\Gamma(\cD^r) \cap \cR^\Gamma(\cN)$ belongs to $H\cdot
\cR^\Gamma(\cN)$. As $H\cdot \cR^\Gamma(\cD ^r)$ is the submodule
$\sum u^{C_i}\cR^\Gamma(\cD ^r)$, the vector-operator $Q$, which
belongs to $\cN$, can be written  $ Q= \sum_{i=1}^p Q_i$ for some
$Q_i \in V^\Gamma_{s-C_i}(\cD ^r)$.

For any $i,j=1,\ldots,p$ we have (see notations \ref{sigmaL}):
$$\ord^{L_i}(Q_j) \leq
L_i(s)-\delta_{ij}.$$
$$\sigma^{L_i}_{L_i(s)}(\sum_j Q_j) = \sigma^{L_i}_{L_i(s)}
(\sum_{j\neq i} Q_j).$$

There exists a family $\ell,\ell_1,\ldots,\ell_p$ of non negative
integer numbers such that $$t^\ell h(Q) = \sum_i t^{\ell_i}
h(Q_i).$$

Let us consider $H_1,\ldots,H_r$ a simultaneous reduced
$L$-standard basis of $h(\cN)$ for any $L$ in $\Gamma$. The
existence of such a simultaneous $L$-standard basis is a
consequence of  the inclusion of $\Gamma$ in a cone of the
analytic standard fan of $h(\cN)$ (see the proof of \cite[Th.
20]{assi-castro-granger-jpaa-2}). By the analytic division theorem
in $\cD[t]$ (see \cite[Th. 7]{assi-castro-granger-jpaa-2}, and
\cite[Th. 4.1.]{granger-oaku}) we can write $$ t^\ell h(Q) =
\sum_{m=1}^r A_m H_m $$ for some $F$--homogeneous elements $A_j
\in \cD[t]$.

For an element $B=\sum_{\alpha \beta} b_{\alpha \beta i} x^\alpha
\partial^\beta e_i$ in $\cD ^r$ we denote by $\cN_V(B)$ its
$V_{\bullet}[\un{n} ]$-Newton diagram. By definition this is the
subset of $\ZZ ^k$ defined by~:
$$\cN_V(B)= \{
(\beta_1-\alpha_1+n^{(i)}_1,\ldots,\beta_k-\alpha_k+n^{(i)}_k)\in
\ZZ^k\,\vert\, b_{\alpha \beta i}\neq 0\}.$$ If no confusion is
posible we shall denote $\beta-\alpha+n^{(i)}=
(\beta_1-\alpha_1+n^{(i)}_1,\ldots,\beta_k-\alpha_k+n^{(i)}_k)$.
We have an analogous definition for the $V_{\bullet}[\un{n}
]$-Newton diagram  of an element in $\cD[t]$.

We have $$\cN_V(t^\ell h(Q))= \cN_V(Q) \subset \bigcup_{i=1}^p
((s-C_i)-\check{\Gamma}).$$ As $H_m$ is $L$-homogeneous for all
$L\in \Gamma$, and as a consequence of the division theorem in
$\cD[t]$, we have $\cN(A_m)+V_{\bullet}[\un{n} ](\exp(H_m))
\subset \bigcup_{i=1}^p ((s-C_i)-\check{\Gamma})$ for
$m=1,\ldots,r$, where $V_{\bullet}[\un{n} ](\alpha,\beta ,i)$ if
by definition  $\beta-\alpha +n^{(i)}$ and $\exp$ stands for the
privileged exponent with respect to a previously fixed well
ordering on $\NN^{2n+1}\times \{ 1\cdots r\}$ (see \cite[pp.
174-175]{granger-oaku} for more details).

As a consequence of the division theorem we also have
$$\sigma^{L_1}_{L_1(s)}(t^\ell h(Q)) = \sum_{j=2}^p
\sigma^{L_1}_{L_1(s)} (t^{\ell_j} h(Q_j)) =
\sum_{m=1}^r\sigma^{L_1}_{L_1(s)} (A_mH_m) = \sum_{m=1}^r
a_m(x,\partial,t) \sigma^{L_1}_{d_{1,m}}(H_m)
$$ where $d_{1,m} = \ord^{L_1}(H_m)$ and $a_m(x,\partial,t)$ is
$L_1$-homogeneous. As the $V_{\bullet}[\un{n} ]$-Newton diagram
of $\sigma^{L_1}_{L_1(s)}(t^\ell h(Q)) = \sum_{j=2}^p
\sigma^{L_1}_{L_1(s)} (t^{\ell_j} h(Q_j))$ is included in
$\bigcup_{j=2}^p ((s-C_j)-\check{\Gamma})$ we have, for any
$m=1,\ldots,r$ the inclusion
$$\cN_V(a_m(x,\partial,t)) +
V_{\bullet}[\un{n} ](\exp(H_m)) \subset \bigcup_{j=2}^p
((s-C_j)-\check{\Gamma}).$$

Let us write $$R_1= t^\ell h(Q) - \sum_{m=1}^r
a_m(x,\partial,t)H_m, \,\, {\mbox{ and }}\,\, Q_1'=(R_1)_{| t=1}.$$ We have :\\
$$\ord^{L_1} (R_1) \leq L_1(s)-1,$$ $$ \ord^{L_j}(R_1) \leq
L_j(s), j=2,\ldots,p,$$ and we can decompose each
$a_m(x,\partial,t)$ into a sum
$$a_m(x,\partial,t) = \sum_{j=2}^p A_{mj}$$ such that
$\cN_V(A_{mj}) + V_{\bullet}[\un{n} ](\exp(H_m)) \subset
((s-C_j)-\check{\Gamma})$ for $m=1,\ldots,r$ and $j=2,\ldots,p.$

Let us write $R_j = \sum_{m=1}^r A_{mj}H_m\in h(\cN)$ and $Q_j'=
(R_j)_{|t=1} \in \cN$ for $j=2,\ldots,p.$

We set $Q_1' = Q - \sum_{j=2} Q_j' \in \cN$. We have $Q_1' \in
V^{\Gamma }_{s-C_1}(\cN)$ by the properties of the orders $
\ord^{L_j}(R_1)$, and for $i=2,\cdots , k$ $Q_i' \in V^{\Gamma
}_{s-C_i}(\cN)$, by the properties of the  $V_{\bullet}[\un{n}
]$-Newton diagram of the $A_{mj}$. Therefore we can write
$$Qu^s=(R_1)_{|t=1}=Q_1'u^{s-C_1}u^{C_1} + \cdots + Q_p'u^{s-C_p}u^{C_p} \in H
\cR^\Gamma(\cN ).$$ That ends the proof.
\end{prueba}


\begin{thebibliography}{11}


\bibitem{Assi-Castro-Granger}
A. Assi., F. J. Castro-Jim\'{e}nez and M. Granger,
\newblock {\rm How to calculate
the slopes of a $D$-module},
\newblock Compositio Math., 104 (1996) 107-123.

\bibitem{assi-castro-granger-jpaa}
A. Assi., F. Castro-Jim\'{e}nez and M. Granger,
\newblock {\rm The Gr\"obner  fan of a $A_n$-module},
\newblock J. Pure Appl. Algebra 150 (1) (2000) 107-123.


\bibitem{assi-castro-granger-jpaa-2}
A. Assi., F. Castro-Jim\'{e}nez and M. Granger,
\newblock {\rm The analytic standard  fan of a $D$-module},
\newblock J. Pure Appl. Algebra 164 (2001) 3-21.

\bibitem{bahloul}
R. Bahloul,
\newblock {\rm D\'emonstration constructive de l'existence de polyn\^omes
de Bernstein-Sato},
\newblock Compositio Math. 41 (2005) 175-191.

\bibitem{bahloul-takayama}
R. Bahloul, N, Takayama,
\newblock {\rm Local Gr\"obner fan : polyedral and computational approach},
\newblock Preprint arXiv:math.AG/0412044, December 2004.

%
%
%
%
%
%
%


\bibitem{granger-oaku} M. Granger, T. Oaku, \newblock{\rm Minimal
filtered free resolutions for analytic D-modules},
\newblock {Journal of pure and applied algebra 191, (2004) 157-180.}


%
%
%
%
%
%
%


%

\bibitem{Herrmann-Orbanz} M. Herrmann and U. Orbanz, \newblock Two notes on flatness.
Manuscripta Math. 40 (1982), no. 1, 109--133.

\bibitem{Sab-Cas}
C. Sabbah and F. Castro,
\newblock {\rm Appendice \`a ``proximit\'e evanescente'' I. La structure
polaire d'un $\cD$--module},
\newblock Compositio Math. 62  (1987) 320-328.

\bibitem{Sabbah}
C. Sabbah, \newblock {\rm Proximit\'e evanescente I. La structure
polaire d'un $\cD$--module}
\newblock Compositio Math. 62  (1987) 283-319.

\end{thebibliography}
\end{document}